\documentclass[graybox]{svmult}

\usepackage{type1cm}        
\usepackage{makeidx}         
\usepackage{graphicx}        
\usepackage{multicol}        
\usepackage[bottom]{footmisc}
 \usepackage{amsmath}

\usepackage{newtxtext}       %
\usepackage{newtxmath}       

\newcommand{\cqfd}{\hfill $\square$}

\newcommand{\R}{\mathbb R}

 \newcommand\revision[1]{{\color{black}{{#1}}}}
 
\makeindex             

\begin{document}

\title*{On the behavior of extreme \revision{$d$-dimensional} spatial quantiles under minimal assumptions}
\author{Davy Paindaveine and Joni Virta}
\institute{Davy Paindaveine \at Universit\'{e} libre de Bruxelles (ECARES and Department of Mathematics) and Universit\'{e} Toulouse~1 Capitole (Toulouse School of Economics), Av. F.D. Roosevelt, 50, CP114/04, 1050, Brussels, Belgium, \email{dpaindav@ulb.ac.be}
\and Joni Virta \at University of Turku (Department of Mathematics and Statistics) and Aalto University School of Science (Department of Mathematics and Systems Analysis), 20014 Turun yliopisto, Finland, \mbox{\email{joni.virta@utu.fi}}}
%
%
\maketitle

\abstract*{}

\abstract{\emph{Spatial} or \emph{geometric} quantiles are the only multivariate quantiles coping with both high-dimensional data and functional data, also in the framework of multi\-ple-output quantile regression. \revision{This work studies spatial quantiles in the finite-dimensional case, where the} spatial quantile~$\mu_{\alpha,u}(P)$ of \revision{the distribution~$P$ taking values in $ \mathbb{R}^d $} is a point in~$\R^d$ indexed by an order~$\alpha\in[0,1)$ and a direction~$u$ in the unit sphere~$\mathcal{S}^{d-1}$ of~$\R^d$---or equivalently by a vector~$\alpha u$ in the open unit ball of~$\R^d$. Recently, \cite{GirStu2017} proved that (i) the extreme quantiles~$\mu_{\alpha,u}(P)$ obtained as~$\alpha\to 1$ exit all compact sets of~$\R^d$ and that (ii) they do so in a direction converging to~$u$. These results help understanding the nature of these quantiles: the first result is particularly striking as it holds even if~$P$ has a bounded support, whereas the second one clarifies the delicate dependence of spatial quantiles on~$u$. However, they were established under assumptions imposing that~$P$ is non-atomic, so that it is unclear whether they hold for empirical probability measures. We improve on this by proving these results under much milder conditions, allowing for the sample case. This prevents using gradient condition arguments, which makes the proofs very challenging. We also weaken the well-known sufficient condition for uniqueness of \revision{finite-dimensional} spatial quantiles.}


\section{Introduction}  
\label{sec:intro}

The problem of defining a satisfactory concept of multivariate quantiles \revision{in $ \mathbb{R}^d $} is a classical one and has generated a huge literature in nonparametric statistics; we refer to \cite{Ser2002C} and the references therein. One of the most famous solutions is given by the \emph{spatial} or \emph{geometric} quantiles introduced in \cite{Cha1996}, which are a particular case of the multivariate M-quantiles from \cite{BreCha1988}; see also \cite{Kol1997}. Spatial quantiles are defined as follows. 

\begin{definition}
Let~$P$ be a probability measure over~$\R^d$. Fix~$\alpha\in[0,1)$ and~$u\in\mathcal{S}^{d-1}$, where~$\mathcal{S}^{d-1}:=\{z\in\R^d:\|z\|^2=z'z=1\}$ is the unit sphere in~$\R^d$. We will say that~$\mu_{\alpha,u}=\mu_{\alpha,u}(P)$ is a \emph{spatial quantile of order~$\alpha$ in direction~$u$ for~$P$} if and only if it minimizes  the objective function 
$$
\mu
\mapsto
O^P_{\alpha,u}(\mu)
:=
\int_{\R^d}
\big\{ 
\| z - \mu \| - \| z \| - \alpha u' \mu
\big\}
\,
dP(z)
$$	
over~$\R^d$ $($the second term in the integrand may look superfluous as it does not depend on~$\mu$, but it actually allows avoiding any moment conditions on~$P)$. 
\end{definition}

Existence and uniqueness of~$\mu_{\alpha,u}$ will be discussed in the next section. It is easy to check \revision{that, for~$d=1$, spatial quantiles reduce to the usual univariate quantiles}. The success of spatial quantiles is partly  explained by their ability to cope with high-dimensional data and even functional data; see, e.g., \cite{Caretal2017}, \cite{Caretal2013}, \cite{ChaCha2014B} and~\cite{ChaCha2014}. These quantiles were also used with much success to conduct multiple-output quantile regression, again also in the framework of functional data analysis; we refer to \cite{ChaLai2013}, \cite{CheGoo2007}, and \cite{ChoCha2019}. \revision{The present work, however, focuses on the finite-dimensional case.} 

In a slightly different perspective, spatial quantiles allow measuring the centrality of any given location in~$\R^d$ with respect to the probability measure~$P$ at hand: if the location~$z$ in~$\R^d$ coincides with the quantile~$\mu_{\alpha,u}$, then a centrality measure for~$z$ is given by its \emph{spatial depth}~$1-\alpha$; see \cite{Gao2003}, \cite{Ser2002A} or \cite{VarZha2000}. This also leads to a spatial concept of multivariate ranks; see, e.g., \cite{Ser2010}. For recent results on spatial depth and spatial ranks, we refer to \cite{Ser2019a,Ser2019b} and to the references therein. The deepest point of~$P$, equivalently its most central quantile, is the quantile~$\mu_0:=\mu_{0,u}$ obtained for~$\alpha=0$ (the dependence on~$u$ of course vanishes at~$\alpha=0$). This is the celebrated \emph{spatial median}, which is one of the earliest robust location functionals; see, e.g.,  \cite{Bro1983} or \cite{Hal1948}. For the other quantiles, the larger~$\alpha$ is, the less central the quantiles~$\mu_{\alpha,u}$ are in each direction~$u$. 

The focus of the present work is on the extreme spatial quantiles that are obtained as~$\alpha$ converges to one. Recently, Girard and Stupfler~\cite{GirStu2017} derived striking results on the behaviour of such extreme spatial quantiles; see also \cite{GirStu2015}. In particular, they showed that, under some assumptions on~$P$ that do not require that~$P$ has a bounded support, these quantiles exit all compact sets of~$\R^d$. Their results, however, require in particular that~$P$ is non-atomic, hence remain silent about empirical distributions~$P_n$ associated with a random sample of size~$n$ from~$P$. Of course, consistency results will imply that the behaviour of sample extreme quantiles will mimic the behaviour of the corresponding population quantiles as~$n$ diverges to infinity; yet for any fixed~$n$, even for large~$n$, there is no guarantee that the results of \cite{GirStu2017} will apply. The goal of the present work is therefore to establish some of these results on extreme spatial quantiles under less stringent assumptions, that will allow for the sample case. Beyond this, we will also weaken the well-known sufficient condition for uniqueness of spatial quantiles. Our results are stated and discussed in Section~\ref{sec:results}, then are proved in Section~\ref{sec:proofs}.



\section{Results}  
\label{sec:results}

We will say that~$P$ is concentrated on a line with direction~$u_*(\in\mathcal{S}^{d-1})$ if and only if there exists~$z_0\in\R^d$ such that~$P[\{z_0+\lambda u_*:\lambda\in\R\}]=1$. Of course, we will say that~$P$ is concentrated on a line if and only if there exists~$u_*\in\mathcal{S}^{d-1}$ such that~$P$ is concentrated on a line with direction~$u_*$. We then have the following existence and uniqueness result.

\begin{theorem}
\label{TheorExistUnicity}
Let~$P$ be a probability measure over~$\R^d$. Fix~$\alpha\in[0,1)$ and~$u\in\mathcal{S}^{d-1}$. Then, (i) $P$ admits a spatial quantile~$\mu_{\alpha,u}$. (ii) If~$P$ is not concentrated on a line, then~$\mu_{\alpha,u}$ is unique. (iii) If~$P$ is not concentrated on a line with direction~$u$, then~$\mu_{\alpha,u}$ is unique for any~$\alpha>0$. \revision{(iv) If~$P$ is concentrated on a line with direction~$u$, say, the line~$\mathcal{L}=\{z_0+\lambda u, \lambda\in\R\}$, then any spatial quantile~$\mu_{\alpha,u}$ belongs to~$\mathcal{L}$; in this case, any such quantile is of the form~$\mu_{\alpha,u}=z_0+\ell_{\alpha}u$, where~$\ell_{\alpha}$ is a spatial quantile of order~$\alpha$ in direction~$1$ for~$P_{z_0,u}$, with~$P_{z_0,u}$ the distribution of~$u'(Z-z_0)$ when~$Z$ has distribution~$P$.}    
\end{theorem}

The existence result in Theorem~\ref{TheorExistUnicity}(i) was established by \cite{Kem1987}, but, since this paper is not easily accessible, we provide our own proof in Section~\ref{sec:proofs}. The uniqueness result in Theorem~\ref{TheorExistUnicity}(ii) is well-known and can be proved by generalizing to an arbitrary quantile the proof for the median in \cite{MilDuc1987}. The result in Theorem~\ref{TheorExistUnicity}(iii) is original and shows that the only case where uniqueness of~$\mu_{\alpha,u}$,~$\alpha>0$, may fail is the one where~$P$ is concentrated on a line with the corresponding direction~$u$. If~$P$ is indeed of this form, then uniqueness may fail exactly as for univariate \revision{(spatial)} quantiles; for instance, if~$P$ is the uniform distribution on~$\{(-2,0),(-1,0),(0,0),(1,0),(2,0)\}$, then any point of the form~$(z,0)$ with~$1\leq z\leq 2$ is a spatial quantile of order~$\alpha=.6$ in direction~$u=(1,0)$ \revision{(recall that the indexing of the classical univariate quantiles differs from the center-outward indexing used for spatial quantiles). Finally, note that, in case~(iii), the spatial quantile~$\mu_{\alpha,u}$ may belong to the line on which~$P$ is concentrated (an example is given below the proof of Lemma~\ref{LemUnicity}).} 

Our main goal is to establish, under very mild conditions, two results  that were recently proved in \cite{GirStu2017} under the assumptions that~$P$ is non-atomic and is not concentrated on a line. The first result states that, as~$\alpha$ converges to one, spatial quantiles with order~$\alpha$ will exit all compact sets in~$\R^d$. Our extension of this result is the following. 

\begin{theorem}
\label{TheorNorm}
Let~$P$ be a probability measure over~$\R^d$. Let~$(\alpha_n)$ be a sequence in~$[0,1)$ that converges to one and let~$(u_n)$ be a sequence in $\mathcal{S}^{d-1}$. Assume that, for any accumulation point~$u_*$ of~$(u_n)$, $P$ is not concentrated on a line with direction~$u_*$ or 
\begin{equation}
	\label{momenttt}
\int_{\R^d} (\|z\|+u_*'z) \,dP(z)
=
\infty
.
\end{equation}
Then, $\|\mu_{\alpha_n,u_n}\|\to\infty$ as~$n\to\infty$ for any sequence of quantiles~$(\mu_{\alpha_n,u_n})$.  
\end{theorem}

Some comments are in order. First, the result does not require that spatial quantiles are unique, which materializes in the fact that the result is stated "for any sequence of quantiles". \revision{Second, the} result allows for distributions that are concentrated on a line, provided that the "moment-type" \revision{Condition~(\ref{momenttt}) is satisfied. Clearly, it is necessary that~$P$ has infinite first-order moments (hence, an unbounded support) for this condition to be satisfied. It is not sufficient, though, as can be seen by considering the limiting behaviour, as~$\alpha\to 1$, of~$\mu_{\alpha,u}$ for a probability measure that would be the distribution of the random vector~$Z=-|\Lambda| u$, where~$\Lambda$ is Cauchy.}
\revision{Third, note that} the result applies as soon as~$P$ is not concentrated on (typically, a few) specific lines, namely those with a direction given by an accumulation point of~$(u_n)$. For instance, if~$u_n=u$ for any~$n$, then the result applies in particular as soon as~$P$ is not concentrated on a line with direction~$u$. But this condition is not even necessary, as the above Cauchy example shows: for instance, in the Cauchy example above, $\|\mu_{\alpha,-u}\|\to\infty$ as~$\alpha\to 1$. Last but not least, Theorem~\ref{TheorNorm} does not require that~$P$ is non-atomic.

\begin{figure}[htbp!] 
\hspace*{-4mm}
\includegraphics[width=1.07\textwidth]{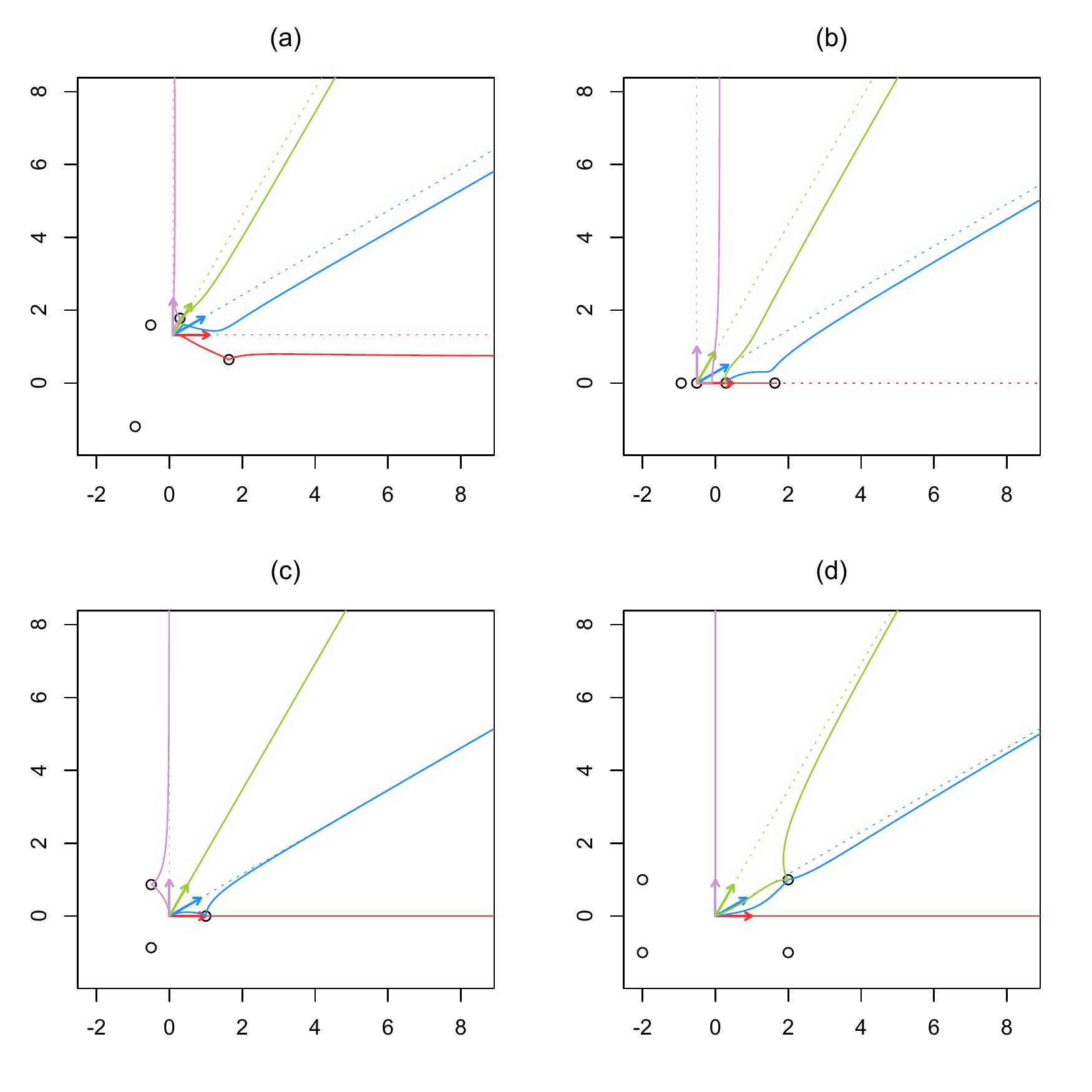} 
\vspace{-4mm}
\caption{For $u=(\cos (\pi j/6),\sin (\pi j/6))$, with~$j=0$ (red), $1$ (blue), $2$ (green) and $3$ (purple), the plots show (linear interpolations of) the spatial quantiles~$\mu_{\alpha_m,u}$, $\alpha_m=.001,.002,\ldots,.999$, in each of the examples~(a)--(d) described in Section~\ref{sec:results}. Dashed lines are showing the halflines with corresponding directions~$u$ originating from the spatial median.}
\label{Fig1}
\end{figure}

We illustrate this result on the basis of the following four examples, in which~$P=P_n$ is the empirical measure associated with a sample~$z_1,\ldots,z_n\in\R^2$. In Example~(a), $n=4$ and the $z_i$'s   were randomly drawn from the \revision{uniform distribution over~$[-2,2]^2$.}
The $z_i$'s in Example~(b) are obtained by projecting those in Example~(a) onto the line~$\{(\lambda,0):\lambda\in\R\}$, whereas those in Example~(c) are~$z_i=(\cos \theta_i,\sin \theta_i)$, $i=1,2,3$, with~$\theta_i=2\pi i/3$, hence are the vertices of an equilateral triangle. Finally, the four~$z_i$'s in Example~(d) are the vertices~$(\pm 2,\pm 1)$ of a rectangle. \revision{These four settings were chosen since they represent point patterns in general position, along a line, on the vertices of a regular polygon, and on the vertices of a stretched regular polygon, respectively.} For each of these examples, Figure~\ref{Fig1} shows the corresponding~$z_i$'s as well as, for four different directions~$u$ (namely, $u=(\cos (\pi j/6),\sin (\pi j/6))$, $j=0,1,2,3$), (linear interpolations of) the spatial quantiles~$\mu_{\alpha_m,u}$, $\alpha_m=.001,.002,\ldots,.999$. The results are perfectly in line with Theorem~\ref{TheorNorm}. Note in particular that, in Example~(b), in which~$P$ is concentrated on the line with direction~$u_*=(1,0)$, the spatial quantiles~$\mu_{\alpha,u}$ exit all compact sets of~$\R^2$ when~$u\neq (\pm) u_*$, as anticipated by Theorem~\ref{TheorNorm}. This fails to happen for~$u=u_*$, which is the only case in Figure~\ref{Fig1} for which our theoretical result remains silent.


The second result from \cite{GirStu2017} we generalize essentially states that the extreme spatial quantiles~$\mu_{\alpha,u}$ are eventually to be found in direction~$u$, which gives a clear interpretation to the direction~$u$ in which quantiles are considered (the directions of non-extreme spatial quantiles do not allow for such a clear interpretation). Our version of this result is the following.

\begin{theorem}
\label{TheorDirection}
Let~$P$ be a probability measure over~$\R^d$. Let~$(\alpha_n)$ be a sequence in~$[0,1)$ that converges to one and let~$(u_n)$ be a sequence in $\mathcal{S}^{d-1}$ that converges to~$u$. Assume that~$P$ is not concentrated on a line with direction~$u$ or that
$$
\int_{\R^d} (\|z\|+u'z) \,dP(z)
=
\infty
$$
Then, 
$
\mu_{\alpha_n,u_n}/\|\mu_{\alpha_n,u_n}\| \to u
$
as~$n\to\infty$ for any sequence of quantiles~$(\mu_{\alpha_n,u_n})$.   
\end{theorem}

The same comments made below Theorem~\ref{TheorNorm} can be repeated here, but for the fact that the sequence~$(u_n$) here may only have one accumulation point, namely its limit~$u$. Again, the result holds for atomic probability measures, which allows us to illustrate the results in Examples~(a)--(d) above. Clearly, Figure~\ref{Fig1} reflects well the conclusion of Theorem~\ref{TheorDirection} in all cases, including those where the probability measure~$P$ is concentrated on a line (again, the case associated with~$u=(1,0)$ in Example~(b) is the only one for which our result remains silent).



\section{Proofs}
\label{sec:proofs}

The proof of Theorem~\ref{TheorExistUnicity} requires \revision{the following three lemmas.}

\begin{lemma}
\label{LemConvexity}
Let~$P$ be a probability measure over~$\R^d$. Fix~$\alpha\in[0,1)$ and~$u\in\mathcal{S}^{d-1}$. Then, (i)~$\mu\mapsto O^P_{\alpha,u}(\mu)$ is convex over~$\R^d$, that is, for $\mu_0,\mu_1\in\R^d$ $(\mu_0\neq \mu_1)$ and~$t\in(0,1)$, one has~$O^P_{\alpha,u}(\mu_t)\leq  (1-t) O^P_{\alpha,u}(\mu_0)+t O^P_{\alpha,u}(\mu_1)$, where we let~$\mu_t:=(1-t) \mu_0+t \mu_1$. (ii) With the same notation, if $P$ is not concentrated on the line containing~$\mu_0$ and~$\mu_1$, then~$O^P_{\alpha,u}(\mu_t)<(1-t) O^P_{\alpha,u}(\mu_0)+t O^P_{\alpha,u}(\mu_1)$.    
\end{lemma}

{\sc Proof of Lemma~\ref{LemConvexity}.}
Fix~$\mu_0,\mu_1\in\R^d$ and $t\in(0,1)$. Then, with~$\mu_t=(1-t) \mu_0+t \mu_1$, we readily have \revision{
\begin{align}\label{forconvex}
\begin{split}
& \| z - \mu_t \| - \| z \| - \alpha u' \mu_t \\
\leq &
(1-t) \{ \| z - \mu_0 \| - \| z \| - \alpha u' \mu_0 \}
+
t \{ \| z - \mu_1 \| - \| z \| - \alpha u' \mu_1 \}.
\end{split}
\end{align}}Part~(i) of the result is then obtained by integrating over~$\R^d$ with respect to~$P$. As for Part~(ii), it follows from the fact that the inequality in~(\ref{forconvex}) is strict for any~$z$ that does not belong to the line containing~$\mu_0$ and~$\mu_1$.
\cqfd
\vspace{3mm}


\revision{
\begin{lemma}
	\label{LemFormerTheorem1i}
	Let~$P$ be a probability measure over~$\R^d$. Fix~$\alpha\in[0,1)$ and~$u\in\mathcal{S}^{d-1}$. Then, $P$ admits a spatial quantile~$\mu_{\alpha,u}$.
\end{lemma}

{\sc Proof of Lemma~\ref{LemFormerTheorem1i}.}
Write~$B_R:=\{z\in\R^d: \|z\|\leq R\}$ and fix~$\lambda>(1+\alpha)/(1-\alpha)$. Pick~$R_0$ large enough so that $P[B_{R_0}]\geq \lambda/(\lambda+1)$. Then,
$$
O^P_{\alpha,u}(\mu)
=
\int_{\R^d}
\big\{ 
\| z- \mu \| 
-
\| z \| 
-
\alpha u'\mu
\big\}
\,
dP(z)
=
O_1(\mu)+O_2(\mu),
$$
where we have
\revision{\begin{align*}
	O_1(\mu)
	:=&
	\int_{B_{R_0}}
	\big\{ 
	\| z- \mu \| 
	-
	\| z \| 
	-
	\alpha  u'\mu
	\big\}
	\,
	dP(z)\\[2mm]
	\geq &
	\int_{B_{R_0}}
	\big\{ 
	\| \mu \|  
	-
	2 \| z \| 
	-
	\alpha \|\mu\|
	\big\}
	\,
	dP(z)\\[2mm]
	\geq &
	\frac{\lambda(1-\alpha)\| \mu \| }{\lambda+1} - 2R_0 
	\end{align*}}
and
\revision{\begin{align*}
	O_2(\mu)
	:=&
	\int_{\R^d\setminus B_{R_0}}
	\big\{ 
	\| z- \mu \| 
	-
	\| z \| 
	-
	\alpha  u'\mu
	\big\}
	\,
	dP(z)\\[2mm]
	\geq&
	\int_{\R^d\setminus B_{R_0}}
	\big\{ 
	- \| \mu \| 
	-
	\alpha \|\mu\|
	\big\}
	\,
	dP(z)\\[2mm]
	\geq &
	- \frac{(1+\alpha)\| \mu \|}{\lambda+1} 
	\cdot
	\end{align*}}
Therefore, for any~$\mu$, we have 
$$
O^P_{\alpha,u}(\mu)
\geq
\frac{\lambda(1-\alpha)-(1+\alpha)}{\lambda+1} 
\| \mu \|
- 2R_0
=:
c_{\lambda,\alpha} \| \mu \| - 2R_0
,
$$
where~$c_{\lambda,\alpha}$ is strictly positive. 
To conclude, pick~$R>0$ so that~$c_{\lambda,\alpha}R-2R_0>O^P_{\alpha,u}(0)$. As a convex function, $\mu \mapsto O^P_{\alpha,u}(\mu)$ is continuous, hence admits a minimum, $\mu_*$ say, in the compact set~$K:=\{\mu\in\R^d: \|\mu\|\leq R\}$. Since any~$\mu\notin K$ is such that
$$
O^P_{\alpha,u}(\mu)
\geq
c_{\lambda,\alpha} R - 2R_0>O^P_{\alpha,u}(0)\geq \min_{\mu\in K} O^P_{\alpha,u}(\mu)
,
$$
we conclude that~$\mu_*$ also minimizes~$\mu\mapsto O^P_{\alpha,u}(\mu)$ over~$\R^d$, which establishes the result.   

\cqfd
\vspace{3mm} 
}


\begin{lemma}
\label{LemUnicity}
Let~$P$ be a probability measure over~$\R^d$ that is concentrated on a line, $\mathcal{L}$ say, with direction~$u_*\in\mathcal{S}^{d-1}$. Fix~$\alpha\in(0,1)$ and~$u\in\mathcal{S}^{d-1}\setminus\{\pm u_*\}$. Then, either~$\mu_{\alpha,u}$ is unique and belongs to~$\mathcal{L}$, or there exists a quantile~$\mu_{\alpha,u}$ that does not belong to~$\mathcal{L}$. 
\end{lemma}

{\sc Proof of Lemma~\ref{LemUnicity}.}
\revision{By Lemma \ref{LemFormerTheorem1i}, there exists at least a quantile~$\mu_{\alpha,u}$.} Trivially, the same proof also shows that~$\mu\mapsto O^P_{\alpha,u}(\mu)$ has a minimizer on~$\mathcal{L}$. Fix then~$\mu_*(\in\mathcal{L})$ arbitrarily such that~$O^P_{\alpha,u}(\mu_*)\leq O^P_{\alpha,u}(\mu)$ for any~$\mu\in\mathcal{L}$.

Let~$Z$ be a random $d$-vector with distribution~$P$. By assumption, $Z=\mu_*+\Lambda u_*$ for some random variable~$\Lambda$, with distribution~$P^\Lambda$ say. For any~$v\in\mathcal{S}^{d-1}$ and any~$h>0$, we then have
\revision{\begin{align*}
\frac{O^P_{\alpha,u}(\mu_*+hv)-O^P_{\alpha,u}(\mu_*)}{h}
&=
- \alpha u'v
+
\int_{\R^d}
\frac{
	\| z - (\mu_*+hv) \| - \| z - \mu_* \| 
}{h}
\,
dP(z) \\
&=
- \alpha u'v
+
\int_{\R}
\frac{
	\| \lambda u_* - hv \| - \| \lambda u_* \| 
}{h}
\,
dP^{\Lambda}(\lambda)
,
\end{align*}}
so that
\revision{$$
\frac{O^P_{\alpha,u}(\mu_*+hv)-O^P_{\alpha,u}(\mu_*)}{h}
-
\big\{
P^{\Lambda}[\{0\}] - s_P u_*'v - \alpha u'v	
\big\}
= 
\int_{\R}
\ell_h(\lambda)
\,
dP^{\Lambda}(\lambda)
,
$$}
where we let 
$
s_P
:=
{\rm E}[{\rm Sign}(\Lambda)]
$
and
$$
\ell_h(\lambda)
:=
\frac{
\| \lambda u_* - hv \| - \| \lambda u_* \| 
}{h}
-
\Big\{
\mathbb{I}[\lambda=0] - {\rm Sign}(\lambda)  (u_*'v) \mathbb{I}[\lambda\neq 0]
\Big\}
.
$$
It is easy to check that, for any~$\lambda\in\R$, the limit of~$\ell_h(\lambda)$ as~$h\to 0$ from above exists and is equal to zero. Moreover, by using the inequality $|\|x\|-\|y\||\leq\|x-y\|$, it is readily seen that the function~$\lambda\mapsto \revision{|\ell_h(\lambda)|}$ is upper-bounded by the function~$\lambda \mapsto 2+|u_*'v|$ that does not depend on~$h$ and is~$P^\Lambda$-integrable. Therefore, Lebesgue's Dominated Convergence Theorem entails that~$\mu\mapsto O^P_{\alpha,u}(\mu)$ admits a directional derivative in direction~$v$ at~$\mu_*$, and that this directional derivative is given by 
\begin{equation}
	\label{dirder}
\frac{\partial O^P_{\alpha,u}}{\partial v}(\mu_*)
=
P^{\Lambda}[\{0\}] - v'(s_P u_* + \alpha u)
.
\end{equation}

Now, using the fact that~$u_*$ and~$u$ are linearly independent and that~$\alpha>0$, one has
$$
m_{\alpha,u}(\mu_*)
:=
\min_{v\in\mathcal{S}^{d-1}}
\frac{\partial O^P_{\alpha,u}}{\partial v}(\mu_*)
=
P^{\Lambda}[\{0\}] - \|s_P u_* + \alpha u\|
,
$$
where the minimum is reached at~$v_0:=(s_P u_* + \alpha u)/\|s_P u_* + \alpha u\|(\neq u_*)$ only. We then consider two cases. 
(i) $m_{\alpha,u}(\mu_*)<0$: then, there exists~$h>0$ such that~$O^P_{\alpha,u}(\mu_*+hv_0)<O^P_{\alpha,u}(\mu_*)$, in which case~$O^P_{\alpha,u}(\mu_*+hv_0)<O^P_{\alpha,u}(\mu)$ for any~$\mu\in\mathcal{L}$, so that any global minimizer of~$\mu\mapsto O^P_{\alpha,u}(\mu)$ does not belong to~$\mathcal{L}$. 
(ii) $m_{\alpha,u}(\mu_*)\geq 0$: then, any directional derivative in~(\ref{dirder}) associated with~$v\in\mathcal{S}^{d-1}\setminus \{v_0\}$ is strictly positive, so that, for any such~$v$, one has~$O^P_{\alpha,u}(\mu_*+hv)>O^P_{\alpha,u}(\mu_*)$  for any~$h$ in an interval of the form~$(0,\varepsilon_v)$. \revision{Pick then, for a fixed~$v\in\mathcal{S}^{d-1}\setminus \{v_0\}$ and the corresponding interval~$(0,\varepsilon_v)$, an arbitrary $ h \in [\varepsilon_v, \infty) $ and any $ h_\varepsilon \in (0,\varepsilon_v)$, and write $ h_\varepsilon = (1 - \lambda) \times 0 + \lambda h $, for $ \lambda := h_\varepsilon/h \in (0, 1) $. The convexity of $ O^P_{\alpha,u} $ (Lemma~\ref{LemConvexity}(i)) entails that
\[ 
\lambda \{O^P_{\alpha,u}(\mu_*+hv) - O^P_{\alpha,u}(\mu_*)\} \geq  O^P_{\alpha,u}(\mu_*+h_\varepsilon v) - O^P_{\alpha,u}(\mu_*) > 0,
\]	
showing that actually~$O^P_{\alpha,u}(\mu_*+hv)>O^P_{\alpha,u}(\mu_*)$ for any~$ h > 0 $.} Continuity of~$\mu\mapsto O^P_{\alpha,u}(\mu)$ (which also follows from convexity) implies that $f(h):=O^P_{\alpha,u}(\mu_*+hv_0)\geq f(0)$ for any~$h>0$ (would there exist~$h>0$ such that~$O^P_{\alpha,u}(\mu_*+hv_0)-O^P_{\alpha,u}(\mu_*)=f(h)-f(0)<0$, then, from continuity, there would exist~$v\in\mathcal{S}^{d-1}\setminus\{v_0\}$ such that~$O^P_{\alpha,u}(\mu_*+hv)-O^P_{\alpha,u}(\mu_*)<0$, a contradiction). It follows that~$\mu_*$ minimizes~$\mu\mapsto O^P_{\alpha,u}(\mu)$ over~$\R^d$. If~$f(h)>f(0)$ for any~$h>0$, then this minimizer is unique, whereas if~$O^P_{\alpha,u}(\mu_*+h_0 v_0)=f(h_0)=f(0)=O^P_{\alpha,u}(\mu_*)$ for some~$h_0>0$, then~$\mu_*+h_0 v_0\notin \mathcal{L}$ also minimizes~$\mu\mapsto O^P_{\alpha,u}(\mu)$ over~$\R^d$. The result follows.   
%
%
\cqfd
\vspace{3mm} 

In the framework of Lemma~\ref{LemUnicity}, it may indeed happen that~$\mu_{\alpha,u}$ is unique and belongs to~$\mathcal{L}$. For instance, if~$P$ is the uniform distribution over~$\{(-1,0),(0,0),(1,0)\}\subset\R^2$, $\alpha\in(0,\frac{1}{3})$ and~$u=(0,1)$, then~$P$ is concentrated on the line~$\mathcal{L}=\{\lambda u_*:\lambda\in\R\}$, with~$u_*=(1,0)$, and~$\mu_{\alpha,u}=(0,0)\in\mathcal{L}$ is the unique order-$\alpha$ quantile in direction~$u$ for~$P$ (this can be checked by proceeding as in the proof of Lemma~\ref{LemUnicity}). 


We can now prove Theorem~\ref{TheorExistUnicity}.
\vspace{2mm}


{\sc Proof of Theorem~\ref{TheorExistUnicity}.}
(i) \revision{The result is an exact restatement of Lemma \ref{LemFormerTheorem1i}.}

(ii) The proof is a straightforward extension of the one in \cite{MilDuc1987}. By contradiction, assume that there exist~$\mu_0$ and~$\mu_1$, with~$\mu_0\neq\mu_1$, such that~$O^P_{\alpha,u}(\mu_0)=O^P_{\alpha,u}(\mu_1)$ is the minimum of~$\mu\mapsto O^P_{\alpha,u}(\mu)$ over~$\R^d$. Since, by assumption,~$P$ is not concentrated on the line containing~$\mu_0$ and~$\mu_1$, Lemma~\ref{LemConvexity}(ii) readily yields that, for any~$t\in(0,1)$,    
$$
O^P_{\alpha,u}((1-t)\mu_0+t\mu_1)<(1-t) O^P_{\alpha,u}(\mu_0) + t O^P_{\alpha,u}(\mu_1) = O^P_{\alpha,u}(\mu_0),
$$   
which contradicts the fact that~$\mu_0$ minimizes~$\mu\mapsto O^P_{\alpha,u}(\mu)$. 

(iii) 
As in the proof of Part~(ii), assume by contradiction that~$\mu\mapsto O^P_{\alpha,u}(\mu)$ has at least two minimizers in~$\R^d$, now with~$\alpha>0$. In view of Part~(ii) of the result, it is enough to consider the case where~$P$ would be concentrated on a line~$\mathcal{L}$ with direction~$u_*(\neq \pm u)$. Lemma~\ref{LemUnicity} thus applies and guarantees that there exists a minimizer of~$\mu\mapsto O^P_{\alpha,u}(\mu)$ that does not belong to~$\mathcal{L}$. Thus, it is possible to pick minimizers~$\mu_0$ and~$\mu_1$ of~$\mu\mapsto O^P_{\alpha,u}(\mu)$, with~$\mu_0\notin\mathcal{L}$ and~$\mu_0\neq\mu_1$. Clearly, $P$ is not concentrated on the line containing~$\mu_0$ and~$\mu_1$ (would it be the case, then~$P$ would be the Dirac measure at the intersection, $\{\mu\}$ say, between~$\mathcal{L}$ and the line containing~$\mu_0$ and~$\mu_1$, hence in particular would be concentrated on the line~$\{\mu+\lambda u:\lambda\in\R\}$ that has direction~$u$, a contradiction). Therefore, Lemma~\ref{LemConvexity}(ii) again yields that, for any~$t\in(0,1)$,    
$$
O^P_{\alpha,u}((1-t)\mu_0+t\mu_1)<(1-t) O^P_{\alpha,u}(\mu_0) + t O^P_{\alpha,u}(\mu_1) = O^P_{\alpha,u}(\mu_0),
$$  
which contradicts the fact that~$\mu_0$ minimizes~$\mu\mapsto  O^P_{\alpha,u}(\mu)$.

\revision{(iv) Assume that~$P$ is concentrated on~$\mathcal{L}=\{z_0+\lambda u, \lambda\in\R\}$. Fix~$\mu\notin\mathcal{L}$. Let us first show that~$\mu$ is not a spatial quantile of order~$\alpha$ in direction~$u$ for~$P$. To do so, write~$Z=\mu_\mathcal{L}+\Lambda u$, where~$\mu_\mathcal{L}$ is the orthogonal projection of~$\mu$ onto~$\mathcal{L}$. Define further~$w:=(\mu_\mathcal{L}-\mu)/c$, with~$c:=\|\mu_\mathcal{L}-\mu\|$. Since~$u'w=0$, we then have 
\begin{eqnarray*}
\lefteqn{
\hspace{-8mm}
\frac{O^P_{\alpha,u}(\mu+hw)-O^P_{\alpha,u}(\mu)}{h}
=
- \alpha u'w
+
\int_{\R^d}
\frac{
	\| z - (\mu+hw) \| - \| z - \mu \| 
}{h}
\,
dP(z) 
}
\\[2mm]
& &
\hspace{10mm}
=
\int_{\R}
\frac{
	\| (\mu_\mathcal{L}+\lambda u) - (\mu+hw) \| - \| (\mu_\mathcal{L}+\lambda u) - \mu \| 
}{h}
\,
dP^{\Lambda}(\lambda)
\\[2mm]
& &
\hspace{10mm}
=
\int_{\R}
\frac{
	\| \lambda u+cw-hw \| - \| \lambda u+cw \| 
}{h}
\,
dP^{\Lambda}(\lambda)
.
\end{eqnarray*}
This yields
$$
\frac{O^P_{\alpha,u}(\mu+hw)-O^P_{\alpha,u}(\mu)}{h}
+
\int_{\R}
\frac{w'(\lambda u+cw)}{\| \lambda u+cw \|}
\,
dP^{\Lambda}(\lambda)
 =
\int_{\R}
g_h(\lambda)
\,
dP^{\Lambda}(\lambda)
,
$$
where 
\begin{eqnarray*}
g_h(\lambda)
&:=&
\frac{
	\| \lambda u+cw-hw \| - \| \lambda u+cw \| 
}{h}
+
\frac{w'(\lambda u+cw)}{\| \lambda u+cw \|}
\\[2mm]
&=&
\frac{
	h^2 - 2 hw'(\lambda u+cw)  
}{h(\| \lambda u+cw-hw \| + \| \lambda u+cw \|)}
+
\frac{w'(\lambda u+cw)}{\| \lambda u+cw \|}
\cdot
\end{eqnarray*}
Clearly, $\lambda\mapsto |g_h(\lambda)|$ is, for~$h\in(0,1)$ say, upper-bounded by the function~$\lambda\mapsto (1/\| \lambda u+cw \|)+3$ that is~$P^\Lambda$-integrable and does not depend on~$h$ (integrability follows from the fact that~$\| \lambda u+cw \|^2=\lambda^2+c^2\geq c^2$). Moreover,  $g_h(\lambda)\to 0$ as~$h\to 0$ for any~$\lambda$. Lebesgue's Dominated Convergence Theorem thus shows that the directional derivative of~$O^P_{\alpha,u}$ at~$\mu$ in direction~$w$ exists and is equal to
$$
\frac{\partial O^P_{\alpha,u}}{\partial w}(\mu)
=
-
\int_{\R}
\frac{w'(\lambda u+cw)}{\| \lambda u+cw \|}
\,
dP^{\Lambda}(\lambda)
=
-
\int_{\R}
\frac{c}{\| \lambda u+cw \|}
\,
dP^{\Lambda}(\lambda)
<
0
.
$$
Therefore, 
$\mu$ is not a spatial quantile of order~$\alpha$ in direction~$u$ for~$P$. 

Consequently, all spatial quantiles of order~$\alpha$ in direction~$u$ for~$P$ belong to~$\mathcal{L}$. These can be characterized as follows. Redefine the random variable~$\Lambda$ through~$Z=z_0+\Lambda u$ (in other words,~$\Lambda=u'(Z-z_0)$). Spatial quantiles are the minimizers of~$\mu\mapsto O^P_{\alpha,u}(\mu)$ over~$\R^d$, which (we just showed it) coincide with the minimizers of the same mapping over~$\mathcal{L}$. These minimizers take the form~$z_0+\ell_{\alpha} u$, where~$\ell_{\alpha}$ minimizes
$$
\lambda\mapsto O^P_{\alpha,u}(z_0+\lambda u)
=
\int_{\R^d}
\{\| z-(z_0+\lambda u) \| - \| z \| -\alpha u'(z_0+\lambda u)\}
\,
dP(z)
$$
$$
=
-\alpha u'z_0
+
\int_{\R}
\{|t-\lambda| - \|z_0+t u\| -\alpha \lambda\}
\,
dP^{\Lambda}(t)
,
$$
or, equivalently, minimizes
$$
\lambda\mapsto 
\int_{\R}
\{|t-\lambda| - |t| -\alpha \lambda\}
\,
dP^{\Lambda}(t)
$$
(note that this last (objective) function, hence also the corresponding minimizers, do not depend on~$u$, which a posteriori justifies the notation~$\ell_{\alpha}$). In other words, $\ell_{\alpha}$ is a spatial quantile of order~$\alpha$ in direction~$1$ for~$P^\Lambda$.  
}
\cqfd
\vspace{3mm}


The proof of Theorem~\ref{TheorNorm} requires both following preliminary results. 

%
%
%
%
%
%

\begin{lemma}
\label{LemContinuity}
Let~$P$ be a probability measure over~$\R^d$. Then, the function
\begin{equation}
\label{objcont}	
(\alpha,u,\mu)
\mapsto 
O^P_{\alpha,u}(\mu)
=
\int_{\R^d}
\big\{ 
\| z - \mu \| - \| z \| - \alpha u' \mu
\big\}
\,
dP(z)
\end{equation}
is continuous over~$[0,1] \times \mathcal{S}^{d-1} \times \R^d$. 
\end{lemma}

{\sc Proof of Lemma~\ref{LemContinuity}.}
Since
\begin{eqnarray*}
	\lefteqn{
\hspace{-3mm} 
	|O^P_{\alpha_2,u_2}(\mu_2)-O^P_{\alpha_1,u_1}(\mu_1)|
}
\\[2mm]
& &
\hspace{3mm} 
\leq
\int_{\R^d}
\big| 
\| z - \mu_2 \| - \| z - \mu_1 \| - (\alpha_2 u_2' \mu_2 - \alpha_1 u_1' \mu_1)
\big|
\,
dP(z)
\\[2mm]
& &
\hspace{3mm} 
\leq
\| \mu_2 - \mu_1 \| 
+
|\alpha_2 u_2' \mu_2 - \alpha_1 u_1' \mu_1|
\\[2mm]
& &
\hspace{3mm} 
\leq
\|\mu_2\| |\alpha_2 - \alpha_1| 
+
 \|\mu_2\| \|u_2-u_1\|
+
(1+\alpha_1) \| \mu_2 - \mu_1 \| 
,
\end{eqnarray*}
the function in~(\ref{objcont}) is Lipschitz over any bounded subset of~$[0,1] \times \mathcal{S}^{d-1} \times \R^d$. The result follows. 
\cqfd
\vspace{3mm}

\begin{lemma}
\label{Lempositionalpha1}
Let~$P$ be a probability measure over~$\R^d$ and fix~$u\in\mathcal{S}^{d-1}$. Assume that~$P$ is not concentrated on a line with direction~$u$ or that
\begin{equation}
	\label{infeq}
\int_{\R^d} (\|z\|+u'z) \,dP(z)
=
\infty
.
\end{equation}
Then the function
$$
\mu 
\mapsto
O^P_{1,u}(\mu)
:=
\int_{\R^d}
\big\{ 
\| z - \mu \| - \| z \| - u' \mu
\big\}
\,
dP(z)
$$
does not have a minimum in~$\R^d$. 
\end{lemma}

{\sc Proof of Lemma~\ref{Lempositionalpha1}.}
Since~$P$ and~$u$ are fixed, we will write~$g(\mu):=O^P_{1,u}(\mu)$ throughout the proof. Letting~$\mu_n:=nu$ (with~$n$ a positive integer), this allows us to write
\begin{eqnarray*}
	g(\mu_n)
&=&
\int_{\R^d}
\big\{ 
\| z - nu \| - (\| z \| + n) 
\big\}
\,
dP(z)
\\[2mm]
&=&
- 2n
\int_{\R^d}
\frac{\|z\|+u'z}{\| z - nu \| + \| z \| + n} 
\,
dP(z)
\\[3mm]
&=&
g_<(\mu_n)+g_\geq(\mu_n)
,
\end{eqnarray*}
where we let
$$
	g_<(\mu_n)
:=
- 2n
\int_{\R^d}
\frac{(\|z\|+u'z) \mathbb{I}[u'z< 0]}{\| z - nu \| + \| z \| + n} 
\,
dP(z)
\
(\leq 0)
$$
and
$$
g_\geq(\mu_n)
:=
- 2n
\int_{\R^d}
\frac{(\|z\|+u'z) \mathbb{I}[u'z\geq 0]}{\| z - nu \| + \| z \| + n} 
\,
dP(z)
\
(\leq 0)
.
$$
 
Now, note that if~(\ref{infeq}) holds, then 
$$
\int_{\R^d} \|z\| \mathbb{I}[u'z\geq 0]\,dP(z)
=
\infty
\ \
\textrm{ or }
\
\int_{\R^d} (\|z\|+u'z) \mathbb{I}[u'z< 0]\,dP(z)
=
\infty
$$
(or both integrals are infinite). This leads to consider three cases. 
\vspace{3mm}
 
Case~(A): $\int_{\R^d} \|z\| \mathbb{I}[u'z\geq 0]\,dP(z)=\infty$. Of course, we have
$$
-g_\geq(\mu_n)
\geq
2n
\int_{\R^d}
\frac{\|z\| \mathbb{I}[u'z\geq 0]}{\| z - nu \| + \| z \| + n} 
\,
dP(z)
.
$$
Since
$
(\| z \| + n)^2
-
\| z - nu \|^2
=
2n\| z \|+2n u'z
\geq 0
$,
we also have
\begin{equation}
	\label{dsd}
-g_\geq(\mu_n)
\geq
\int_{\R^d}
\frac{n\|z\| \mathbb{I}[u'z\geq 0]}{\| z \| + n} 
\,
dP(z)
=:
\int_{\R^d}
h_n(z) 
\,
dP(z)
.
\end{equation}
Since $h_n(z)\leq h_{n+1}(z)$ for any~$z$ and the pointwise limit of~$h_n$ is the function~$h$ defined by~$h(z):=\|z\|\mathbb{I}[u'z\geq 0]$, the Monotone Convergence Theorem yields 
$$
\int_{\R^d}
h_n(z) 
\,
dP(z)
\to 
\int_{\R^d}
h(z)
\,
dP(z)
=
\infty
,
$$
which, jointly with~(\ref{dsd}), establishes that~$g_\geq(\mu_n)\to -\infty$. Since~$g(\mu_n)\leq g_\geq(\mu_n)$, we conclude that~$g(\mu_n)\to -\infty$, so that~$g$ does not have a minimum in Case~(A).
\vspace{3mm}
    

Case~(B): $\int_{\R^d}
(\|z\|+u'z) \mathbb{I}[u'z< 0]
\,
dP(z)
=
\infty$. Using the Monotone Convergence Theorem as in Case~(A) readily provides that
\begin{eqnarray*}
-
g_<(\mu_n)
&=&
2n
\int_{\R^d}
\frac{(\|z\|+u'z) \mathbb{I}[u'z< 0]}{\| z - nu \| + \| z \| + n} 
\,
dP(z)
\\[4mm]
&=&
2
\int_{\R^d}
\frac{(\|z\|+u'z) \mathbb{I}[u'z< 0]}{\sqrt{\frac{1}{n^2}\|z\|^2 + 1 + \frac{2}{n}|u'z|} + \frac{1}{n}\| z \| + 1} 
\,
dP(z)
\end{eqnarray*}
converges to
$$
\int_{\R^d}
(\|z\|+u'z) \mathbb{I}[u'z< 0]
\,
dP(z)
=
\infty
$$
as~$n\to\infty$. Since~$g(\mu_n)\leq g_<(\mu_n)$, this yields~$g(\mu_n)\to -\infty$. It follows that~$g$ does not have a minimum in Case~(B).  
\vspace{3mm}


Case~(C): $\int_{\R^d} \|z\| \mathbb{I}[u'z\geq 0]\,dP(z)<\infty$ and~$\int_{\R^d}
(\|z\|+u'z) \mathbb{I}[u'z< 0]
\,
dP(z)
< \infty$. Using the finiteness of the first and second integrals, Lebesgue's \revision{Dominated Convergence Theorem} readily yields
$$
g_\geq(\mu_n)
\to 
- 
\int_{\R^d}
(\|z\|+u'z) \mathbb{I}[u'z\geq 0] 
\,
dP(z)
$$
and
$$
g_<(\mu_n)
\to
- 
\int_{\R^d}
(\|z\|+u'z)
 \mathbb{I}[u'z< 0] 
\,
dP(z)
,
$$
respectively. Therefore,
$$
	g(\mu_n)
	=
	g_<(\mu_n)+g_\geq(\mu_n)
\to 
- 
\int_{\R^d}
(\|z\|+u'z) 
\,
dP(z)
=:
i_u^P
.
$$
In Case~(C), $P$ is not concentrated on a line with direction~$u$ by assumption, which implies that, for any~$\mu\in\R^d$,
$$
g(\mu)
-
i_u^P
=
\int_{\R^d}
\big\{ 
\| z - \mu \| + u' (z-\mu)
\big\}
\,
dP(z)
>
0
.
$$
This shows that the function $g$ does not have a minimum in Case~(C) either. The result is thus proved.
\cqfd
\vspace{3mm}

Theorem~\ref{TheorNorm} then follows from Lemmas~\ref{LemContinuity}--\ref{Lempositionalpha1} in the same way as Theorem~2.1(i) in \cite{GirStu2017} (but for the fact that we are considering distributions that do not ensure uniqueness of quantiles). We still report the proof for the sake of completeness. 
\vspace{3mm}
 
{\sc Proof of Theorem~\ref{TheorNorm}.}
Ad absurdum, assume that there exists a sequence of quantiles~$(\mu_{\alpha_n,u_n})$ such that~$\|\mu_{\alpha_n,u_n}\|$ does not diverge to infinity. Then,  $(\mu_{\alpha_n,u_n},u_n)$ has a subsequence that is bounded, hence from compactness, possesses a further subsequence, $(\mu_{\alpha_{n_\ell},u_{n_\ell}},u_{n_\ell})$ say,  that converges in~$\R^d\times \mathcal{S}^{d-1}$, to~$(\mu_\infty,u_\infty)$, say. By construction, $u_\infty$ is an accumulation point of the sequence~$(u_n)$. For any~$\ell$, we have
$$
O_{\alpha_{n_\ell},u_{n_\ell}}^P(\mu_{\alpha_{n_\ell},u_{n_\ell}})
\leq
O_{\alpha_{n_\ell},u_{n_\ell}}^P(\mu)
$$
for any~$\mu\in\R^d$. In view of Lemma~\ref{LemContinuity}, taking limits as~$\ell\to\infty$ then provides
$$
O_{1,u_\infty}^P(\mu_\infty)
\leq
O_{1,u_\infty}^P(\mu)
$$
for any~$\mu\in\R^d$. Since this contradicts Lemma~\ref{Lempositionalpha1}, the result is proved. 
\cqfd
\vspace{3mm}


The proof of Theorem~\ref{TheorDirection} requires the following lemma.

\begin{lemma}
	\label{LemDirection}
Let~$P$ be a probability measure over~$\R^d$ and fix~$m\in(0,2)$. Then, 
$$
t_P(r)
:=
\int_{\R^d}
\frac{
\|z\| 
}
{
\sqrt{ (\|z\|-r)^2+m r\|z\|}
}
\,
dP(z)
\to 0
$$
as~$r\to\infty$. 
\end{lemma}

{\sc Proof of Lemma~\ref{LemDirection}.}
Fix~$\delta>0$. For any~$r>0$, let~$Y_r:=\|Z\|/r$, where~$Z$ is a random $d$-vector with distribution~$P$. Then, with~$h:=m\delta^2/4$,  
\begin{eqnarray*}
t_P(r)
&=&
{\rm E}\Bigg[
\frac{
\|Z\| 
}
{
\sqrt{ (\|Z\|-r)^2+m r\|Z\|}
}
\Bigg]
=
{\rm E}\Bigg[
\frac{
Y_r
}
{ 
\sqrt{ (Y_r-1)^2+m Y_r}
}
\Bigg]  
\\[3mm] 
&=&
{\rm E}\Bigg[
\frac{ 
Y_r \mathbb{I}[Y_r\leq h]
}   
{
\sqrt{ (Y_r-1)^2+m Y_r}  
} 
\Bigg]
+ 
{\rm E}\Bigg[ 
\frac{  
Y_r \mathbb{I}[Y_r> h]
}
{
\sqrt{ (Y_r-1)^2+m Y_r}
}
\Bigg]
.
\end{eqnarray*}
Since~$y/\sqrt{(y-1)^2+my}\leq 2/\sqrt{m(4-m)}$ for any~$y\geq 0$, this provides
\begin{eqnarray*}
t_P(r)
& \leq &
{\rm E}\Bigg[
\frac{
\sqrt{Y_r} \mathbb{I}[Y_r\leq h]
}
{
\sqrt{m}
}
\Bigg]
+
\frac{2}{\sqrt{m(4-m)}}
P[Y_r> h]
\\[2mm]
& \leq &
\frac{\delta}{2}
+
\frac{2}{\sqrt{m(4-m)}}
P[\|Z\|> rh]
<
\delta
,
\end{eqnarray*}
for~$r$ large enough. 
\cqfd
\vspace{3mm}


{\sc Proof of Theorem~\ref{TheorDirection}.}
In this proof, we use the notation
$$
\mathcal{S}^{\rm in}_{u,c}
:= 
\mathcal{S}^{d-1}\cap \{ z\in\R^d : u'z \geq 1-c \}
$$ 
and
$$
\mathcal{S}^{\rm out}_{u,c}
:= 
\mathcal{S}^{d-1}\cap \{ z\in\R^d : u'z \leq 1-c \}
.
$$
Ad absurdum, assume that there exists a sequence of quantiles~($\mu_{\alpha_n,u_n}$) such that $(w_n:=\mu_{\alpha_n,u_n}/\|\mu_{\alpha_n,u_n}\|)$ does not converge to~$u$. Thus, there exists~$\varepsilon>0$ such that 
$
w_n
\in
\mathcal{S}^{\rm out}_{u,\varepsilon}
$ 
for infinitely many~$n$. Upon extraction of a subsequence, we may assume that~$w_n$ belongs to~$\mathcal{S}^{\rm out}_{u,\varepsilon}$ for any~$n$. By assumption, we may, still upon extraction of a subsequence, assume that $u_n\in\mathcal{S}^{\rm in}_{u,\varepsilon/2}$ for any~$n$. Assume for a moment that there exist~$R>0$ and~$\eta\in(0,1)$ such that
\begin{equation}
\label{tos}
O^P_{\alpha,v}(rw)
> 
O^P_{\alpha,v}(rv)
\end{equation}
for any~$\alpha\in[\eta,1)$, $r\geq R$, $v\in\mathcal{S}^{\rm in}_{u,\varepsilon/2}$ and~$w\in\mathcal{S}^{\rm out}_{u,\varepsilon}$. Pick then~$n$ large enough to have~$\alpha_n\geq \eta$ and~$\|\mu_{\alpha_n,u_n}\|\geq R$ (existence follows from Theorem~\ref{TheorNorm}). By definition, this implies that
$$
O^P_{\alpha_n,u_n}(\|\mu_{\alpha_n,u_n}\| w_n)
=
O^P_{\alpha_n,u_n}(\mu_{\alpha_n,u_n})
\leq 
O^P_{\alpha_n,u_n}(\|\mu_{\alpha_n,u_n}\| u_n)
,
$$
which contradicts~(\ref{tos}). 
 
Therefore, it is sufficient to prove~(\ref{tos}). To do so, fix~$v\in\mathcal{S}^{\rm in}_{u,\varepsilon/2}$, $w\in\mathcal{S}^{\rm out}_{u,\varepsilon}$ and~$\eta\in(0,1)$ (we show that~(\ref{tos}) holds, actually, not just for some~$\eta\in(0,1)$ but for any~$\eta\in(0,1)$). Note that one has
$
\sqrt{2(1-v'w)} = \|v-w\| \geq u'(v-w)=u'v-u'w\geq (1-\varepsilon/2)-(1-\varepsilon)=\varepsilon/2
$
so that
$
2(1-v'w)\geq \varepsilon^2/4,
$
hence
$$
v'w\leq 1-\frac{\varepsilon^2}{8}
\cdot
$$
Write then
\begin{eqnarray*}
\lefteqn{
\hspace{-15mm} 
O^P_{\alpha,v}(rw)
-
O^P_{\alpha,v}(rv)
=
\int_{\R^d}
\big\{ 
\| z - rw \| -\| z - rv \|  - \alpha (r v'w-r)
\big\}
\,
dP(z)
}
\\[2mm]
& & 
\hspace{-3mm} 
=
r \alpha (1- v'w)
+
\int_{\R^d}
\frac{
\| z - rw \|^2 - \| z - rv \|^2 
}
{
\| z - rw \| + \| z - rv \| 
}
\,
dP(z)
\\[2mm]
& & 
\hspace{-3mm} 
\geq
\frac{r \eta \varepsilon^2}{8}
+
\int_{\R^d}
\frac{
2r(v-w)'z
}
{
\| z - rv \| + \| z - rw \| 
}
\,
dP(z)
\\[2mm]
& & 
\hspace{-3mm} 
\geq
r 
\Bigg[
\frac{\eta \varepsilon^2}{8}
-4
\int_{\R^d}
\frac{
\|z\| 
}
{
\| z - rv \|+\| z - rw \|
}
\,
dP(z)
\Bigg]
.
\end{eqnarray*}

Now, using the fact that~$\|v+w\|^2=2(1+v'w)\leq 2(2-\varepsilon^2/8)$, we obtain 
\begin{eqnarray*}
\hspace{2mm} 
	\lefteqn{
\{\| z - rv \|
+
\| z - rw \|\}^2
\geq
\| z - rv \|^2
+
\| z - rw \|^2
}
\\[2mm]
& & 
\hspace{-4mm} 
=
2\|z\|^2+2r^2-2r(v+w)'z
\geq
2\|z\|^2+2r^2-2\sqrt{2(2-\varepsilon^2/8)} r\|z\|
\\[2mm]
& & 
\hspace{-4mm} 
=
2\{ (\|z\|-r)^2+\sqrt{2}(\sqrt{2}-\sqrt{2-\varepsilon^2/8})r\|z\|\}
=:
2\{ (\|z\|-r)^2+m_\varepsilon r\|z\|\}
,
\end{eqnarray*}
which provides
$$
O^P_{\alpha,v}(rw)
-
O^P_{\alpha,v}(rv)
\geq
r 
\Bigg[
\frac{\eta \varepsilon^2}{8}
-2\sqrt{2}
\int_{\R^d}
\frac{
\|z\| 
}
{
\sqrt{(\|z\|-r)^2+m_\varepsilon r\|z\|}
}
\,
dP(z)
\Bigg]
.
$$
Since~$m_\varepsilon\in(0,2)$, Lemma~\ref{LemDirection} guarantees that there exists~$R>0$, not depending on the choice of $ v, w, \eta $ and $ \alpha $, such that for any~$r\geq R$, 
$
O^P_{\alpha,v}(rw)
-
O^P_{\alpha,v}(rv)
\geq
r 
\eta \varepsilon^2/16
>
0
$. This proves~(\ref{tos}), hence the result. 
\cqfd

\begin{acknowledgement}
Davy Paindaveine's research is supported by a research fellowship from the Francqui Foundation and by the Program of Concerted Research Actions (ARC) of the Universit\'{e} libre de Bruxelles. The research of Joni Virta was supported by the Academy of Finland (grant 321883).	
\end{acknowledgement}


\bibliography{Paper}           
\bibliographystyle{spmpsci} 
\vspace{3mm}


\end{document}